\documentclass[12pt]{article}

\usepackage{amsmath}
\usepackage{rawfonts} 

\input prepictex
\input pictex
\input postpictex

\newcounter{mathitem}
\newenvironment{mathitem}
  {\begin{list}{$(\roman{mathitem})$}{
   \setcounter{mathitem}{0}
   \usecounter{mathitem}
   \setlength{\topsep}{0pt plus 2pt minus 0pt}
   \setlength{\parskip}{0pt plus 2pt minus 0pt}
   \setlength{\partopsep}{0pt plus 2pt minus 0pt}
   \setlength{\parsep}{0pt plus 2pt minus 0pt}
   \setlength{\leftmargin}{35pt}
   \setlength{\itemsep}{0pt plus 2pt minus 0pt}}}
 {\end{list}}

\def\vejde#1{\unskip\nobreak\hfill\penalty50\hskip1em\hbox{}\nobreak\hfill\hbox{#1}}
\newcommand{\mezera}{\vspace{3.34mm}}
\newcommand{\carka}{,\penalty0\relax}
\newcommand{\zlom}{\penalty0\relax}

\newcounter{citac}[section]
\def\thecitac{\thesection.\arabic{citac}}

\newenvironment{Theorem}{\par\mezera\noindent%
\refstepcounter{citac}{\bf Theorem \thecitac.}\quad\bgroup\sl }
{\egroup\par\mezera\endtrivlist}%

\newenvironment{Theoremcite}[1]{\par\mezera\noindent%
\refstepcounter{citac}{\bf Theorem \thecitac{} {\cite{#1}}.}\quad\bgroup\sl }
{\egroup\par\mezera\endtrivlist}%

\newenvironment{Lemma}{\par\mezera\noindent%
\refstepcounter{citac}{\bf Lemma \thecitac.}\quad\bgroup\sl }
{\egroup\par\mezera\endtrivlist}%

\newenvironment{Corollary}{\par\mezera\noindent%
\refstepcounter{citac}{\bf Corollary \thecitac.}\quad\bgroup\sl}
{\egroup\par\mezera\endtrivlist}%

\newenvironment{Proof}{\par\mezera\noindent%
{\bf Proof.}\quad\bgroup}
{\egroup\vejde{\rule[0pt]{2.5mm}{2.5mm}}\par\mezera\endtrivlist}%

\newenvironment{Conjecture}{\par\mezera\noindent%
\refstepcounter{citac}{\bf Conjecture \thecitac.}\quad\bgroup\sl }
{\egroup\par\mezera\endtrivlist}%

\title{Hamiltonian Cycles in the Square of a Graph}

\author{Jan Ekstein\thanks{supported by the grant 1M0545 of the Czech Ministry of                 
                           Education.}\\
\small Department of Mathematics\\[-0.8ex]
\small University of West Bohemia, Pilsen, Czech Republic\\
\small \texttt{ekstein@kma.zcu.cz}}
      
\begin{document}
\maketitle

\begin{abstract}
We show that under certain conditions the square of the graph obtained by identifying
a vertex in two graphs with hamiltonian square is also hamiltonian. Using this 
result, we prove necessary and sufficient conditions for hamiltonicity of the square
of a connected graph such that every vertex of degree at least three in a~block graph 
corresponds to a cut vertex and any two these vertices are at distance at least four.

\end{abstract}

\section{Introduction and notation}
The graphs considered in this paper are undirected and simple. If $G$ is a~graph, we
denote by $V(G)$ the vertex set of $G$, by $E(G)$ the edge set of $G$. For $x \in 
V(G)$, $d_{G}(x)$ denotes the \emph{degree} of $x$ and $N_{G}(x)$ denotes the
\emph{neighbourhood} of $x$. For $x, y \in V(G)$, $\mbox{dist}_G(x,y)$ denotes the 
\emph{distance} between $x, y$. For $A\subseteq V(G)$, $\langle A\rangle$ denotes the 
subgraph of $G$ induced by $A$.

The \emph{$k$-star} is a tree on $k + 1$ vertices with one vertex of degree $k$, 
called the center, and the others of degree 1, $k = 0, 1, 2, ...$. The graph 
$S(K_{1, 3})$ is the graph $K_{1, 3}$ in which each edge is subdivided once. Given 
sets $A, B$ of vertices, we call $P = x_{0}, ..., x_{k}$ an \emph{$(A, B)$-path} if
$V(P) \cap A = \{x_{0}\}$ and $V(P) \cap B = \{x_{k}\}$, we write $(a, B)$-path
rather than $(\{a\}, B)$-path. For a~graph $G$ we define $V_{i}(G) = \{v \in V(G): 
d(v) = i\}$ and $W(G) = V(G) \setminus V_{2}(G)$. A~\emph{branch} in $G$ is a 
nontrivial path whose ends are in $W(G)$ and whose internal vertices, if any, are of 
degree 2 in $G$.

The \emph{square} of $G$, denoted $G^{2}$, is the graph with the vertex set $V(G)$ in 
which two vertices are adjacent if their distance in $G$ is one or two. We say that
two graphs are \emph{homeomorphic} if they can be turned into isomorphic graphs by
finite number of edge-subdivisions. Let $G'$ be an subgraph of $G$. We say that $G'$ 
is \emph{maximal} with respect to a given graph property if $G'$ itself has the 
property but no graph $G' + A$ does, for any nonempty subset $A \subseteq E(G) 
\setminus E(G')$.

A connected graph that has no cut vertices is called a \emph{block}. A \emph{block} 
of a graph is a~subgraph that is a block and is maximal with respect to this 
property. The \emph{degree of a~block} $B$ of a graph $G$, denoted by $d(B)$, is the
number of cut vertices of $G$ belonging to $V(B)$. A block of degree 1 is called an
\emph{endblock} of $G$, otherwise it is a \emph{non-end block}. A block is said to be 
\emph{acyclic} if it is isomorphic to one edge, otherwise we say it is \emph{cyclic}. 
The \emph{block graph} of a graph $G$ is the graph $\mbox{Bl}(G)$ such that the
vertices of $\mbox{Bl}(G)$ are the blocks and cut vertices of $G$, and two vertices
are adjacent in $\mbox{Bl}(G)$ if one of them is a~block of~$G$ and the second one is 
its vertex.

Let $G_{1}$, $G_{2}$ be connected graphs, $x \notin V(G_{1})~\bigcup~V(G_{2})$,
$V(G_{1})~\bigcap~V(G_{2})\zlom = \emptyset$, and let $x_{i} \in V(G_{i})$, $i = 1, 2$. 
Then the graph $G$ with vertex set $V(G) = (V(G_{1})~\setminus~\{x_{1}\})~\bigcup 
\zlom (V(G_{2})~\setminus~\{x_{2}\})~\bigcup~\{x\}$ and with edge set $E(G) = E(G_{1}
- x_{1})~\bigcup~E(G_{2} - x_{2})~\bigcup \zlom \{ux |~u \in V(G_{1}), ux_{1} \in
E(G_{1})\}~\bigcup~$$ $$\{vx |~v \in V(G_{2}), vx_{2} \in E(G_{2})\}$ is called the
\emph{connection of the graphs $G_{1}$, $G_{2}$ over the vertices $x_{1}$, $x_{2}$},
 denoted $G = G_{1}[x_{1} = x_{2}]G_{2}$.

Let $G$ be a connected graph such that $G^{2}$ is hamiltonian and let $x \in V(G)$. 
We say that
     \begin{mathitem}
      \item[a)] the vertex $x$ is \emph{of type 1} if there exists a hamiltonian
                 cycle $C$ of $G^{2}$ such that both edges of $C$
                 incident with $x$ are in $G$,
      \item[b)] the vertex $x$ is \emph{of type 2} if $x$ is not of type 1 and
                 there exists a hamiltonian cycle~$C$ of $G^{2}$ such that
                 exactly one edge of $C$ incident with $x$ is in $G$,
      \item[c)] the vertex $x$ is \emph{of type 3} if $x$ is not of type 1 or 2 and
                 there exists a~hamiltonian cycle $C$ of $G^{2}$ such that
                 for some two vertices $u, v \in N_{G}(x)$ is $uv \in E(C)$,
      \item[d)] the vertex $x$ is \emph{of type 4} if $x$ is not of type 1 or 2 or 3.
      \end{mathitem}
We denote $V_{[i]}(G) = \{x \in V(G)~|~x \mbox{ is of type }i\}, i = 1, 2, 3, 4$.

\section{The connection of graphs}
Let us first mention the following result by Fleischner \cite{Fleis} that will be 
used many times in proofs.

\begin{Theoremcite}{Fleis}
  \label{Fleischner}
Let $y$ and $z$ be arbitrarily chosen vertices of a 2-connected graph~$G$. Then 
$G^{2}$ contains a hamiltonian cycle $C$ such that the edges of $C$ incident with $y$ 
are in~$G$ and at least one of edges of $C$ incident with $z$ is in $G$. If $y$ and 
$z$ are adjacent in $G$, then these are three different edges.
\end{Theoremcite}

It is easy to see that Theorem \ref{Fleischner} implies that the square of 
a 2-connected graph is hamiltonian.

The following result shows that, under certain conditions, the square of the 
connection of two graphs with hamiltonian square is also hamiltonian.

\begin{Theorem}
 \label{ConnectTheorem}
 Let $G_{1}$, $G_{2}$ be connected graphs such that $(G_{1})^{2}$, $(G_{2})^{2}$ are 
 hamiltonian, let $x_{i} \in V(G_{i})$, $i = 1, 2$. If
   \begin{mathitem}
     \item[I)] $G = G_{1}[x_{1}=x_{2}]G_{2}$ and $x_{i} \in V_{[1]}               
                                      (G_{i})~\bigcup~V_{[2]}(G_{i})$, $i = 1,2$, or
     \item[II)] $G = G_{1}[x_{1}=x_{2}]K_{2}$, $x_{1} \in V_{[1]}
                               (G_{1})~\bigcup~V_{[2]}(G_{1})~\bigcup~V_{[3]}(G_{1})$
                                and $V(K_{2}) = \{x_{2}, u\}$ or
     \item[III)] $G = G_{1}[x_{1}=x_{2}]G_{2}$, $x_{1} \in V_{[3]}(G_{1})$ and $x_{2} 
                                                       \in V_{[1]}(G_{2})$,
   \end{mathitem}
 then $G^{2}$ is hamiltonian.

 Moreover under the assumptions of I),
    \begin{mathitem}
     \item[a)] if $x_{i} \in V_{[1]}(G_{i})$, $i = 1, 2$,
                 then $x = x_{1} = x_{2} \in V_{[1]}(G)$;
     \item[b)] if $x_{1} \in V_{[1]}(G_{1})$ and $x_{2} \in V_{[2]}(G_{2})$,
               then $x = x_{1} = x_{2} \in V_{[2]}(G)$;
     \item[c)] if $G_{2}$ is 2-connected and $x_{1} \in V_{[1]}(G_{1}) \bigcup 
                                                               V_{[2]}(G_{1})$, then
               $v \in V_{[1]}(G)$ for any \zlom $v~\in~V(G_{2})$, $v \neq x_{2}$;
     \item[d)] if $x_{i} \in V_{[2]}(G_{i})$, $i = 1, 2$,
               then $x = x_{1} = x_{2} \notin V_{[1]}(G) \bigcup V_{[2]}(G)$.
    \end{mathitem}

 Moreover under the assumptions of II),
    \begin{mathitem}
    \item[a)] if $x_{1} \in V_{[1]}(G_{1})$, then $x = x_{1} = x_{2} \in V_{[1]}(G)$;
    \item[b)] if $x_{1} \in V_{[2]}(G_{1})$, then $x = x_{1} = x_{2} \in V_{[2]}(G)$;
    \item[c)] if $x_{1} \in V_{[1]}(G_{1}) \bigcup V_{[2]}(G_{1})$, then $u \in 
                                                                         V_{[2]}(G)$.
    \end{mathitem}
\end{Theorem}
\newpage
\begin{Proof}
 \begin{mathitem}
 \item[I)] Let $x = x_{1} = x_{2}$ and let $C_{1}$, $C_{2}$ be hamiltonian cycles in 
  $(G_{1})^{2}$, $(G_{2})^{2}$ such that $a_{1}x, b_{1}x \in E(G)$ for $a_{1} \in 
  N_{C_{1}}(x)$, $b_{1} \in N_{C_{2}}(x)$, respectively.  Let $a_{2} \in N_{C_{1}}
  (x)$, $a_{1} \neq a_{2}$, and $b_{2} \in N_{C_{2}}(x)$, $b_{1} \neq b_{2}$, and let 
  $P_{a_{1}a_{2}} = C_{1} - x_{1}$ and $P_{b_{2}b_{1}} = C_{2} - x_{2}$. Then 
  $P_{a_{1}a_{2}}$, $P_{b_{2}b_{1}}$ are hamiltonian paths in $(G_{1} - x)^{2}$, 
  $(G_{2} - x)^{2}$, respectively, and the cycle 
  $C = a_{1}P_{a_{1}a_{2}}a_{2}xb_{2}P_{b_{2}b_{1}}b_{1}a_{1}$ is a hamiltonian cycle 
  in $G^{2}$.
 \begin{mathitem}
   \item[a)] If moreover $x_{i} \in V_{[1]}(G_{i})$, $i = 1, 2$, then we can assume 
             that $a_{2}x, b_{2}x \in E(G)$ and therefore $x \in V_{[1]}(G)$.
   \item[b)] If moreover $x_{1} \in V_{[1]}(G_{1})$ and $x_{2} \in V_{[2]}(G_{2})$, 
             then we can assume that $a_{2}x \in E(G)$ and it is obvious that there 
             is no $C_{2}$ such that $b_{2}x \in E(G)$, therefore $x \in V_{[2]}(G)$.
   \item[c)] If moreover $G_{2}$ is 2-connected, then for any $v \in V(G_{2})$, 
             $v \neq x_{2}$, we can assume by Theorem~\ref{Fleischner} without loss 
             of generality that $c_{1}v, c_{2}v \in E(G_{2})$, $c_{1}, c_{2} \in 
             N_{C_{2}}(v)$, $c_{1} \neq c_{2}$. If $v \neq b_{1}$, then $c_{1}v, 
             c_{2}v \in E(C)$, therefore $v~\in~V_{[1]}(G)$. If $v = b_{1}$, then we 
             have $vx_{2} \in V(G_{2})$ and by Theorem~\ref{Fleischner} $x_{2}b_{2} 
             \in E(G)$. Then $\widetilde{C} = 
             a_{1}P_{a_{1}a_{2}}a_{2}xb_{1}P_{b_{1}b_{2}}b_{2}a_{1}$ is also 
             a hamiltonian cycle in~$G^{2}$ and moreover the edges of $C_{2}$ 
             incident with $v = b_{1}$ are in $\widetilde{C}$. Therefore             
             $v \in V_{[1]}(G)$.
   \item[d)] If moreover $x_{i} \in V_{[2]}(G_{i})$, $i = 1, 2$, then it is obvious 
             that there are no $C_{1}$, $C_{2}$ such that $a_{2}x \in E(G)$ or
            $b_{2}x \in E(G)$ and therefore $x \notin V_{[1]}(G) \bigcup V_{[2]}(G)$.
 \end{mathitem}
 \item[II)] \emph{Case 1: $x_{1} \in V_{[3]}(G_{1})$.}\\
   Let $x = x_{1} = x_{2}$ and let $C_{1}$ be a hamiltonian cycle in $(G_{1})^{2}$
   such that $yw \in E(C_{1})$ for some $y, w \in N_{G_{1}}(x)$. Let $P_{yw} = C_{1} 
   - yw$. Then $P_{yw}$ is a hamiltonian path in $(G_{1})^{2}$ and the cycle $C = 
   yP_{yw}wuy$ is a~hamiltonian cycle in $G^{2}$.

   \emph{Case 2: $x_{1} \in V_{[1]}(G_{1}) \bigcup V_{[2]}(G_{1})$.}\\
      Let $x = x_{1} = x_{2}$ and let $C_{1}$ be a hamiltonian cycle in $(G_{1})^{2}$
      such that $yx \in E(G)$ for $y \in N_{C_{1}}(x)$. Let $z \in N_{C_{1}}(x)$, 
      $z \neq y$, and let $P_{zy} = C_{1} - x_{1}$. Then $P_{yz}$ is a~hamiltonian
      path in $(G_{1} - x)^{2}$ and the cycle $C = zP_{zy}yuxz$ is a hamiltonian
      cycle in $G^{2}$.
 \begin{mathitem}
   \item[a)] If moreover $x_{1} \in V_{[1]}(G_{1})$, then we can assume that 
             $xz \in E(G)$ and therefore $x \in V_{[1]}(G)$.
   \item[b)] If moreover $x_{1} \in V_{[2]}(G_{1})$, then it is obvious that there is 
             no $C_{1}$ such that $xz \in E(G)$ and therefore $x \in V_{[2]}(G)$.
   \item[c)] Since $ux \in E(G)$ and $N_{G}(u) = \{x\}$, it is obvious that there is 
             no cycle~$\widetilde{C}$ in $G^{2}$ such that both edges of         
             $\widetilde{C}$ incident with $u$ are in $G$ and therefore 
             $u \in V_{[2]}(G)$.
 \end{mathitem}
 \item[III)] Let $x = x_{1} = x_{2}$, let $C_{1}$ be a hamiltonian cycle in 
   $(G_{1})^{2}$ such that $yw \in E(C_{1})$ for some $y, w \in N_{G_{1}}(x)$ and let 
   $C_{2}$ be a hamiltonian cycle in $(G_{2})^{2}$ such that $ax, bx \in E(G)$ for
   $a, b \in N_{C_{2}}(x)$, $a \neq b$. Let $P_{yw} = C_{1} - yw$ and 
   $P_{ab} = C_{2} - x_{2}$. Then $P_{yw}$, $P_{ab}$ are hamiltonian paths in
   $(G_{1})^{2}$, $(G_{2} - x)^{2}$, respectively, and the cycle $C =
   wP_{wy}yaP_{ab}bw$ is a~hamiltonian cycle in $G^{2}$.
 \end{mathitem}
\end{Proof}

\section{The hamiltonian square of a graph}
This work is motivated by the following result due to El Kadi Abderrezzak, Flandrin 
and Ryj\'{a}\v{c}ek \cite{ElKadiAbdFlaRyj}.

 \begin{Theoremcite} {ElKadiAbdFlaRyj}
  \label{AbdFlaRyj}
  If $G$ is a connected graph such that every induced $S(K_{1,3})$ has at least three 
  edges in a block of degree at most 2, then $G^{2}$ is hamiltonian.
 \end{Theoremcite}

 The following result, originally by Thomassen \cite{Thom}, is a corollary of Theorem  
 \ref{AbdFlaRyj}.

 \begin{Theoremcite} {Thom}
 \label{Thomassen}
  If the block graph of $G$ is a path, then $G^{2}$ is hamiltonian.
 \end{Theoremcite}

Before the presentation of main results we first give the following slight 
strengthening of Theorem \ref{Thomassen} which will be needed in our proofs.

\begin{Theorem}
 \label{ZobThomassen}
 Let $G$ be a graph such that its block graph is a path and let $u_{1}$, $u_{2}$
 be arbitrary vertices which are not cut vertices and are contained in different
 endblocks of~$G$.\\
 Then $G^{2}$ contains a hamiltonian cycle $C$ such that, for $i = 1, 2$,
  \begin{mathitem}
     \item[$\bullet$] if $u_{i}$ is contained in a cyclic block, then both edges of 
                      $C$ incident with $u_{i}$ are in $G$, and   
     \item[$\bullet$] if $u_{i}$ is contained in an acyclic block, then exactly one 
                      edge of $C$ incident with $u_{i}$ is in $G$.
  \end{mathitem}
\end{Theorem}

\begin{Proof}
 If $G$ is a path of length at least 2, then the theorem is obvious. Thus, suppose  
 that $G$ contains at least one cyclic block $B_{1}$ and let $k$ denote the number of 
 blocks of $G$.
 
 \medskip
 
 We prove the theorem by induction on $k$.\\
 \emph{1.} Let $k = 2$, let $B_{2}$ be the second block of $G$, let $x = V(G)$ be the 
           (only) cut vertex of $G$ and let $u_{1}$, $u_{2}$ be arbitrary vertices 
           such that $u_{1} \in V(B_{1})$, $u_{2} \in V(B_{2})$ and $u_{1} \neq x$, 
           $u_{2} \neq x$. The graph $(B_{1})^{2}$ contains a hamiltonian cycle 
           $C_{1}$ such that the edges of $C_{1}$ incident with $u_{1}$ are in 
           $B_{1}$ and at least one of edges of $C_{1}$ incident with $x$ is in 
           $B_{1}$. If $u_{1}$ and $x$ are adjacent in $G$, then these are three 
           different edges by Theorem \ref{Fleischner}. Then we can assume that 
           $x \in V_{[1]}(B_{1}) \bigcup V_{[2]}(B_{1})$ and $G = B_{1}
           [x_{1}=x_{2}]B_{2}$, where $x_{1}$ and $x_{2}$ is the copy of $x$ in 
           $B_{1}$ and $B_{2}$, respectively.
 \begin{mathitem}
     \item[a)] If $B_{2}$ is cyclic, then  the graph $G^{2}$ contains a hamiltonian 
               cycle $C$ such that both edges of $C$ incident with $u_{2}$ are in $G$ 
               by Theorem \ref{ConnectTheorem} Ic) and it is obvious that we can find 
               $C$ such that also both edges of $C$ incident with $u_{1}$ are in $G$.
     \item[b)] If $B_{2} = K_{2} = x_{2}u_{2}$, then the graph $G^{2}$ contains 
               a hamiltonian cycle $C$ such that exactly one edge of $C$ incident 
               with $u_{2}$ is in $G$ by Theorem \ref{ConnectTheorem} IIc) and it is 
               obvious that we can find $C$ such that both edges of $C$ incident with 
               $u_{1}$ are in $G$.
 \end{mathitem}
 
 \medskip
 
 \emph{2.} Suppose the assertion is true for each graph with at most $k$ blocks, let 
           $G$ be a graph with $k + 1$ blocks such that its block graph is a path and 
           let $u_{1}$, $u_{2}$ be arbitrary vertices which are not cut vertices and 
           are contained in different endblocks of $G$, $k \geq 2$.

Let $B_{k+1}$ be the endblock of $G$ containing $u_{2}$. We denote 
$\widetilde{G} = G - V(B_{k+1} - x)$, where $x \in V(B_{k+1})$ is a cut vertex of 
$G$. Then $G = \widetilde{G}[x_{1}=x_{2}]B_{k+1}$, where $x_{1}$ and $x_{2}$ is the
copy of $x$ in $\widetilde{G}$ and $B_{k+1}$, respectively, and we can assume by the 
induction hypothesis that $\widetilde{G}$ contains a hamiltonian cycle $C_{1}$ such 
that if $u_{1}$, $x_{1}$ is contained in a cyclic block, then both edges of $C_{1}$
incident with $u_{1}$, $x_{1}$ are in $\widetilde{G}$, and if $u_{1}$, $x_{1}$ is 
contained in an acyclic block, then exactly one edge of $C_{1}$ incident with 
$u_{1}$, $x_{1}$ is in $\widetilde{G}$, respectively. Then we can assume that 
$x_{1} \in V_{[1]}(\widetilde{G}) \bigcup V_{[2]}(\widetilde{G})$.
 \begin{mathitem}
     \item[a)] If $B_{k+1}$ is cyclic, $u_{2} \in V(B_{k+1})$ and $u_{2} \neq x_{2}$, 
       then the graph $G^{2}$ contains a hamiltonian cycle $C$ such that both edges 
       of $C$ incident with $u_{2}$ are in $G$ by Theorem~\ref{ConnectTheorem}~Ic) 
       and it is obvious that we can find $C$ such that if $u_{1}$ is contained in 
       a cyclic block, then both edges of $C$ incident with $u_{1}$ are in $G$, and
       if $u_{1}$ is contained in an acyclic block, then exactly one edge of $C$
       incident with $u_{1}$ is in $G$.
     \item[b)] If $B_{k+1} = K_{2} = x_{2}u_{2}$, then the graph $G^{2}$ contains 
       a hamiltonian cycle $C$ such that exactly one edge of $C$ incident with             
       $u_{2}$ is in $G$ by Theorem \ref{ConnectTheorem} IIc) and it is obvious that 
       we can find $C$ such that if $u_{1}$ is contained in a cyclic block, then both 
       edges of~$C$ incident with $u_{1}$ are in $G$, and if $u_{1}$ is contained in 
       an acyclic block, then exactly one edge of $C$ incident $u_{1}$ is in $G$.
 \end{mathitem}
\end{Proof}

\begin{figure}[ht]
  $$\beginpicture
   \setcoordinatesystem units <1mm,1mm>
   \setplotarea x from -30 to 30, y from 5 to 5
    \put{$\bullet$} at 0 0
    \put{$v_{1}$} at 0 -4
    \put{$\bullet$} at 15 -7.5
    \put{$v_{9}$} at 15 -11.5
    \put{$\bullet$} at 15 7.5
    \put{$v_{8}$} at 19 7.5
    \put{$\bullet$} at -15 7.5
    \put{$v_{4}$} at -19 7.5
    \put{$\bullet$} at -15 -7.5
    \put{$v_{2}$} at -15 -11.5
    \put{$\bullet$} at  -7.5 15
    \put{$v_{5}$} at -11.5 15
    \put{$\bullet$} at  7.5 15
    \put{$v_{7}$} at 11.5 15
    \setlinear
    \plot  0 0 15 -7.5 15 7.5 0 0 /
    \plot  0 0 -15 7.5 -15 -7.5 0 0 /
    \plot  0 0 -7.5 15 7.5 15 0 0 /
    \put{$\bullet$} at 30 -7.5
    \put{$v_{10}$} at 30 -11.5
    \plot  15 -7.5 30 -7.5 /
    \put{$\bullet$} at -30 -7.5
    \put{$v_{3}$} at -30 -11.5
    \plot  -15 -7.5 -30 -7.5 /
    \put{$\bullet$} at -7.5 30
    \put{$v_{6}$} at -7.5 34
    \plot -7.5 15 -7.5 30 /
  \endpicture$$
  \caption{}
  \label{Vzor}
 \end{figure}
 
We consider the graph in Figure \ref{Vzor}. It is easy to see that the cycle $C = 
v_{1} ,v_{2}, v_{3} \carka v_{4},v_{5}, v_{6}, v_{7}, v_{8}, v_{9}, v_{10}, v_{1}$ is 
a hamiltonian cycle in $G^{2}$ but the induced subgraph $H = 
\langle{\{{v_{1},v_{2},v_{3},v_{5},v_{6},v_{9},v_{10}}}\}\rangle$ is isomorphic to 
$S(K_{1,3})$ and does not have at least three edges in a block of degree at most 2. 
This example shows that the assumptions in Theorem \ref{AbdFlaRyj} are sufficient but 
not necessary. We looked for other conditions implying that the square of a graph is 
hamiltonian.

\section{Main result}
Let $V_{\geq 3}(G) = \{x \in V(G)| d_G(x) \geq 3\}$ and, for $x \in V(G)$, $t_{G}(x)$ 
denotes the number of acyclic non-end blocks of $G$ containing $x$. First of all we 
prove the following lemma.

 \begin{Lemma}
   \label{Lemma}
    Let $G$ be a connected graph with exactly one vertex in $V_{\geq 3}(Bl(G))$ 
    corresponding to a cut vertex $a$ of $G$. If $a$ is contained in at most two 
    acyclic non-end blocks of $G$, then $G^{2}$ contains a hamiltonian cycle $C$ such 
    that if $t_{G}(a) = 0$, then both edges of $C$ incident with $a$ are in $G$, if 
    $t_{G}(a) = 1$, then exactly one edge of $C$ incident with $a$ is in $G$,
    if $t_{G}(a) = 2$, then no edge of $C$ incident with $a$ is in $G$.
   \end{Lemma}

 \begin{Proof}
  Let $r \geq 0$, $s \geq 0$ and $t = t_{G}(a) \geq 0$ denote the number of cyclic 
  blocks, acyclic endblocks and acyclic non-end blocks of $G$ containing $a$, 
  respectively, and choose the notation such that if $r > 0$, then $B_{1},..., B_{r}$ 
  are all cyclic blocks, if $s > 0$, then $B_{r+1},..., B_{r+s}$ are all acyclic 
  endblocks, and if $t > 0$, then $B_{r+s+1},..., B_{r+s+t}$ are all acyclic non-end 
  blocks of $G$ containing the vertex $a$.

  By the assumption, $t \leq 2$ and $r + s + t \geq 3$, hence $r + s > 0$.

  \emph{Case 1: $r = 0$.}

  \noindent
    If $t = 0$, then $G$ is a star and the assertion is obvious. Let $t \geq 1$.

    \emph{Subcase 1.1: $s = 1$.} Then necessarily $t = 2$. Let $B_{1} = au$ and let 
          $b_{2}$, $b_{3}$ be the branch of $\mbox{Bl}(G)$ containing the vertex 
          corresponding to $B_{2}$, $B_{3}$ and denote $H_{2}$, $H_{3}$ the subgraph
          corresponding to $b_{2}$, $b_{3}$, respectively. For $i = 2, 3$, 
          $\mbox{Bl}(H_{i})$ is a path and therefore $(H_{i})^{2}$ is hamiltonian and 
          $a \in V_{[2]}(H_{i})$ by Theorem \ref{ZobThomassen}. If $G_{1} = H_{2}
          [x_{1}=x_{2}]B_{1}$, where $x_{1}$ and $x_{2}$ is the copy of $a$ in $H_{2}$ 
          and $B_{1}$, respectively, then $(G_{1})^{2}$ is hamiltonian and $a \in 
          V_{[2]}(G_{1})$ by Theorem \ref{ConnectTheorem} IIb). Moreover $G = G_{1}
          [y_{1}=y_{2}]H_{3}$,  where $y_{1}$ and $y_{2}$ is the copy of~$a$ in 
          $G_{1}$ and $H_{3}$, respectively, and $G^{2}$ contains  hamiltonian cycle 
          $C$ such that no edge of $C$ incident with $a$ is in $G$ by Theorem 
          \ref{ConnectTheorem} Id).

    \emph{Subcase 1.2: $s \geq 2$.} Then necessarily $1 \leq t \leq 2$. For 
          $i = s + 1, s + 2$, let $H_{i}$ be the same subgraphs as in Subcase 1.1 and 
          let $\widetilde{G} = G - V(H_{s+1} - a) - V(H_{s+2} - a)$. It is obvious 
          that $\widetilde{G}$ is a star and therefore $(\widetilde{G})^{2}$ is 
          hamiltonian and $a \in V_{[1]}(\widetilde{G})$. If $t = 1$, then 
          $G = \widetilde{G}[x_{1}=x_{2}]H_{s+1}$, where $x_{1}$ and $x_{2}$ is the 
          copy of $a$ in $\widetilde{G}$ and $H_{s+1}$, respectively, and $G^{2}$ 
          contains  hamiltonian cycle $C$ such that exactly one edge of $C$ incident 
          with $a$ is in $G$ by Theorem \ref{ConnectTheorem} Ib).
          Let $t = 2$. If $G_{1} = \widetilde{G}[x_{1}=x_{2}]H_{s+1}$, where $x_{1}$ 
          and $x_{2}$ is the copy of $a$ in $\widetilde{G}$ and $H_{s+1}$, 
          respectively, then $(G_{1})^{2}$ is hamiltonian and $a \in V_{[2]}(G_{1})$ 
          by Theorem~\ref{ConnectTheorem}~Ib). Moreover 
          $G = G_{1}[y_{1}=y_{2}]H_{s+2}$, where $y_{1}$ and $y_{2}$ is the copy of 
          $a$ in $G_{1}$ and $H_{s+2}$, respectively, and $G^{2}$
          contains  hamiltonian cycle $C$ such that no edge of $C$ incident with $a$ 
          is in $G$ by Theorem \ref{ConnectTheorem} Id).
     \newpage
     
    \emph{Case 2: $r \geq 1$.}

    \noindent
    For $i = 1, 2, ..., r$, let $b_{i}$ be the branch of $\mbox{Bl}(G)$ containing
    the vertex corresponding to~$B_{i}$. We denote $H_{i}$ the subgraph corresponding 
    to $b_{i}$. The block graph $\mbox{Bl}(H_{i})$ is a path and therefore
    $(H_{i})^{2}$ is hamiltonian and $a \in V_{[1]}(H_{i}) $ either by Theorem 
    \ref{ZobThomassen} or by Theorem~\ref{Fleischner}. Let $b_{r+s+1}$, $b_{r+s+2}$ be 
    the branch of $\mbox{Bl}(G)$ containing the vertex corresponding to $B_{r+s+1}$, 
    $B_{r+s+2}$ and denote $H_{r+s+1}$, $H_{r+s+2}$ the subgraph corresponding to 
    $b_{r+s+1}$, $b_{r+s+2}$, res\-pectively. For $j = r + s + 1, r + s + 2$, 
    $\mbox{Bl}(H_{j})$ is a path and therefore $(H_{j})^{2}$ is hamiltonian and 
    $a \in V_{[2]}(H_{j})$ by Theorem \ref{ZobThomassen}.
    
    Let $G_{1} = G - V(H_{r+s+1} - a) - V(H_{r+s+2} - a)$ and let $\ell$ denote the 
    number of branches of $\mbox{Bl}(G_{1})$. We show that $(G_{1})^{2}$ is 
    hamiltonian and $a \in V_{[1]}(G_{1})$.
    
    We proceed by induction on $\ell$.

    \noindent
    \emph{1.} For $\ell = 1$ obviously $G_{1} = H_{1}$ and the assertion is true.

    \noindent
    \emph{2.} Suppose the assertion is true for each graph such that its block graph 
       contains at most $\ell$ branches, let $G_{1}$ be a graph without acyclic 
       non-end blocks such that its block graph contains $\ell + 1$ branches and with 
       exactly one vertex in $V_{\geq 3}(Bl(G_{1}))$ corresponding to a~cut vertex $a$ 
       of $G_{1}$.

       If $\widetilde{G_{1}} = G_{1} - V(H_{1} - a)$, then 
       $G_{1} = H_{1}[x_{1}=x_{2}]\widetilde{G_{1}}$, where $x_{1}$ and $x_{2}$ is the 
       copy of $a$ in $H_{1}$ and $\widetilde{G_{1}}$, respectively. If 
       $\widetilde{G_{1}} = B_{r+1}$, then $(G_{1})^{2}$ is hamiltonian and 
       $a \in V_{[1]}(G_{1})$ by Theorem \ref{ConnectTheorem} IIa). Otherwise
       $\widetilde{G_{1}}$ is hamiltonian and $a = x_{1} \in V_{[1]}
       (\widetilde{G_{1}})$ by the induction hypothesis. Then $(G_{1})^{2}$ is 
       hamiltonian and $a \in V_{[1]}(G_{1})$ by Theorem \ref{ConnectTheorem} Ia).

    If $G = G_{1}$, then it is obvious that $G^{2}$ contains a hamiltonian cycle $C$ 
    such that both edges of $C$ incident with $a$ are in $G$. If $t = 1$, then 
    $G = G_{1}[y_{1}=y_{2}]H_{r+s+1}$, where $y_{1}$ and $y_{2}$ is the copy of $a$ in 
    $G_{1}$ and $H_{r+s+1}$, respectively, and $G^{2}$ contains hamiltonian cycle~$C$ 
    such that exactly one edge of $C$ incident with $a$ is in $G$ by Theorem 
    \ref{ConnectTheorem} Ib). Let $t = 2$. If $G_{2} = G_{1}[y_{1}=y_{2}]H_{r+s+1}$, 
    where $y_{1}$ and $y_{2}$ is the copy of $a$ in $G_{1}$ and $H_{r+s+1}$, 
    respectively, then $(G_{2})^{2}$ is hamiltonian and $a \in V_{[2]}(G_{2})$ by 
    Theorem \ref{ConnectTheorem} Ib). Moreover $G = G_{2}[z_{1}=z_{2}]H_{r+s+2}$, 
    where $z_{1}$ and $z_{2}$ is the copy of $a$ in $G_{2}$ and $H_{r+s+2}$, 
    respectively, and $G^{2}$ contains  hamiltonian cycle $C$ such that no edge
    of $C$ incident with $a$ is in $G$ by Theorem \ref{ConnectTheorem} Id).
  \end{Proof}

Now we prove our main result.

\newpage

\begin{Theorem}
 \label{MainTheorem}
  Let $G$ be a connected graph with at least three vertices such that
  \begin{mathitem}
    \item[i)] every vertex $x \in V_{\geq 3}(\mbox{Bl}(G))$ corresponds to a cut 
              vertex of $G$, and
    \item[ii)] for any two vertices $x, y \in V_{\geq 3}(\mbox{Bl}(G))$   
               it holds that $\mbox{dist}_{\mbox{Bl}(G)}(x, y) \geq 4$.
  \end{mathitem}
  Then $G^{2}$ is hamiltonian if and only if every cut vertex of $G$ is contained in 
  at most two acyclic non-end blocks of $G$.
\end{Theorem}

\begin{Proof}
 I) First suppose that we have a vertex $a \in V(G)$ contained in at least three  
    acyclic non-end blocks of $G$. We show that the graph $G^{2}$ is not hamiltonian. 
    Let, to the contrary, $C$ be a hamiltonian cycle in $G^{2}$. For $i = 1, 2, 3$, 
    let $B_{i} = a_{i}a$ denote three acyclic non-end blocks of $G$ and $B_{i+3}$ 
    a~block of $G$ adjacent to the block $B_{i}$ such that $a \notin V(B_{i+3})$. Then 
    necessarily there is a vertex $c_{i} \in N_{B_{i+3}}(a_{i})$ such that 
    $c_{i}a \in E(C)$, for $i = 1, 2, 3$. From this it follows that $d_{C}(a) \geq 3$, 
    contradicting the fact that $C$ is a cycle.

 \medskip

 II) Suppose that every cut vertex of $G$ is contained in at most two acyclic 
     \zlom non-end blocks of $G$. We show that $G^{2}$ is hamiltonian.

    If $G$ is a cyclic block, then $G^{2}$ is hamiltonian by Theorem \ref{Fleischner},
    and if $\mbox{Bl}(G)$ is a path, then $G^{2}$ is hamiltonian by Theorem 
    \ref{Thomassen}.

    Now suppose that $\mbox{Bl}(G)$ contains at least one vertex of degree at least 
    three cor\-responding to a cut vertex of $G$. For $i = 1, 2, ..., k$, let $b_{i}$ 
    be a vertex of $\mbox{Bl}(G)$ in~$V_{\geq 3}(\mbox{Bl}(G))$, let $a_{i}$ be the 
    vertex of $G$ corresponding to $b_{i}$ and we choose the notation such that \zlom
    $\mbox{dist}_{Bl(G)}(b_1, b_k)$ is maximum and the (unique) path in $\mbox{Bl}(G)$
    joining $b_1$ and $b_2$ has no \zlom interior vertices in 
    $V_{\geq 3}(\mbox{Bl}(G))$. Let $t_{G}(a_{i}) \geq 0$ denote the number of acyclic 
    non-end blocks of $G$ containing $a_{i}$.

   We prove the following statement.

   \medskip

   \emph{Under the assumptions of Theorem \ref{MainTheorem} the graph $G^{2}$ contains
   a hamiltonian cycle $C$ such that if $t_{G}(a_{i}) = 0$, then both edges of $C$ 
   incident with $a_{i}$ are in $G$, if $t_{G}(a_{i}) = 1$, then exactly one edge of 
   $C$ incident with $a_{i}$ is in $G$, if $t_{G}(a_{i}) = 2$, then no edge of $C$ 
   incident with $a_{i}$ is in $G$, $i = 1, 2, ..., k$.}

   \medskip

   We proceed by induction on $k$.

   For $k = 1$ the assertion is given by Lemma \ref{Lemma}.

   Suppose the assertion is true for each graph $G'$ such that its block graph 
   $\mbox{Bl}(G')$ has at most $k - 1$ vertices in $V_{\geq 3}(\mbox{Bl}(G'))$ 
   corresponding to cut vertices of $G'$ and these are at distance at least four in 
   $\mbox{Bl}(G')$, and let $G$ be a graph such that its block graph $\mbox{Bl}(G)$
   is a~tree with $k$ vertices in $V_{\geq 3}(\mbox{Bl}(G))$ corresponding to cut 
   vertices of $G$ and any two vertices of $\mbox{Bl}(G)$ in $V_{\geq 3}(\mbox{Bl}
   (G))$ are at distance at least four in $\mbox{Bl}(G)$, $k \geq 2$.

   By the notation of $a_{1}$, let $H$ be the unique subgraph of graph $G$ 
   corresponding to the $(b_{1}, b_{2})$-path in $\mbox{Bl}(G)$. Let 
   $\widetilde{G} = G - V(H - \{a_{1}, a_{2}\})$. We denote the components of 
   $\widetilde{G}$ by $G_{1}$, $G_{2}$ such that $a_{1} \in V(G_{1})$ and $a_{2} \in 
   V(G_{2})$.

   If $d_{\mbox{Bl}(G_{1})}(a_{1}) \geq 3$ and $d_{\mbox{Bl}(G_{2})}(a_{2}) \geq 3$, 
   then, by the induction hypothesis, $(G_{1})^{2}$, $(G_{2})^{2}$ contains 
   a hamiltonian cycle $C_{1}$, $C_{2}$ such that if $t_{G_{1}}(a_{1}) = 0$, 
   $t_{G_{2}}(a_{i}) = 0$, then both edges of $C_{1}$, $C_{2}$ incident with $a_{1}$, 
   $a_{i}$ are in $G_{1}$, $G_{2}$, if $t_{G_{1}}(a_{1}) = 1$, $t_{G_{2}}(a_{i}) = 1$, 
   then exactly one edge of $C_{1}$, $C_{2}$ incident with $a_{1}$, $a_{i}$
   is in $G_{1}$, $G_{2}$, if $t_{G_{1}}(a_{1}) = 2$, $t_{G_{2}}(a_{i}) = 2$, then no 
   edge of $C_{1}$, $C_{2}$ incident with $a_{1}$, $a_{i}$ is in $G_{1}$, $G_{2}$, 
   $i = 2, 3, ..., k$, respectively.

   In the case $d_{\mbox{Bl}(G_{1})}(a_{1}) = 2$, set $K_{2} = vu$ (where 
   $v, u \notin V(G)$), and $\widehat{G_{1}} = G_{1}[a_{1}=v]K_{2}$. Then
   $(\widehat{G_{1}})^{2}$ contains a hamiltonian cycle with the required properties 
   by the induction hypothesis and it is obvious that also $(G_{1})^{2}$. We proceed 
   in the case $d_{\mbox{Bl}(G_{2})}(a_{2}) = 2$ similarly.

   Then by the assumption that any two vertices of $\mbox{Bl}(G)$ in $V_{\geq 3}
   (\mbox{Bl}(G))$ are at distance at least four in $\mbox{Bl}(G)$ and by Theorem 
   \ref{ZobThomassen}, the graph $H^{2}$ contains a hamiltonian cycle $C_{H}$ such 
   that, for $j = 1, 2$, if $a_{j}$ is contained in a cyclic block, then both edges of 
   $C_{H}$ incident with $a_{j}$ are in $H$, and if $a_{j}$ is contained in an acyclic 
   block, then exactly one edge of $C_{H}$ incident with $a_{j}$ is in $H$.

   \medskip

   \emph{Case 1: $t_{G_{2}}(a_{2}) \in \{0, 1\}$.}

  \noindent
  Let $\widetilde{G_{2}} = G_{2}[x_{1}=x_{2}]H$, where $x_{1}$ and $x_{2}$ is the copy
  of $a_{2}$ in $G_{2}$ and $H$, respectively. Then $(\widetilde{G_{2}})^{2}$ contains
  a hamiltonian cycle $\widetilde{C_{2}}$ with the required properties by Theorem 
  \ref{ConnectTheorem} either Ia) or Ib) or Id) (using $C_{H}$ and $C_{2}$). Moreover 
  it is obvious that if $a_{1}$ is contained in a cyclic block of $\widetilde{G_{2}}$, 
  then both edges of $\widetilde{C_{2}}$ incident with $a_{1}$ are in
  $\widetilde{G_{2}}$, and if $a_{1}$ is contained in an acyclic block of 
  $\widetilde{G_{2}}$, then exactly one edge of $\widetilde{C_{2}}$ incident with 
  $a_{1}$ is in~$\widetilde{G_{2}}$.
  
  \medskip
  
  \begin{mathitem}
  \item[a)] If $t_{G_{1}}(a_{1}) \in \{0, 1\}$, then $G = G_{1}
       [y_{1}=y_{2}]\widetilde{G_{2}}$, where $y_{1}$ and $y_{2}$ is the copy of 
       $a_{1}$ in $G_{1}$ and $\widetilde{G_{2}}$, respectively, and $G^{2}$ 
       contains a hamiltonian cycle $C$ such that if $t_{G}(a_{i}) = 0$,
       then both edges of $C$ incident with $a_{i}$ are in $G$, if $t_{G}(a_{i}) = 1$, 
       then exactly one edge of $C$ incident with $a_{i}$ is in $G$, if 
       $t_{G}(a_{i}) = 2$, then no edge of $C$ incident with $a_{i}$ is in $G$, 
       $i = 1, 2, ..., k$, by Theorem \ref{ConnectTheorem} either Ia) or Ib) or Id) 
       (using $\widetilde{C_{2}}$ and $C_{1}$).
       
  \item[b)] Let $t_{G_{1}}(a_{1}) = 2$. Let $B_{1}$ be an acyclic non-end block of 
       $G_{1}$ containing the vertex $a_{1}$, let $F$ be the subgraph of $G_{1}$
       corresponding to the maximal connected subgraph of $\mbox{Bl}(G_{1})$ 
       containing the vertex corresponding to $B_{1}$ and not containing the 
       vertex~$b_{1}$.

       Let $F_{1} = G_{1} - V(F - a_{1})$. By the induction hypothesis, $G_{1}$ 
       contains a hamiltonian cycle $C_{1}$ such that no edge of $C_{1}$ incident with 
       $a_{1}$ is in $G_{1}$ and we can divide $C_{1}$ into hamiltonian cycles 
       $C_{1a}$ in $F_{1}$ and $C_{1b}$ in $F$ such that exactly one edge of $C_{1a}$ 
       incident with $a_{1}$ is in $F_{1}$ and exactly one edge of $C_{1b}$ incident 
       with $a_{1}$ is in $F$.

       Necessarily both edges of $\widetilde{C_{2}}$ incident with $a_{1}$ are in
       $\widetilde{G_{2}}$ (otherwise $t_{G}(a_{1}) = 3$, a~contradiction). Set
       $\widetilde{G_{1}} = F_{1}[y_{1}=y_{2}]\widetilde{G_{2}}$, where $y_{1}$ and 
       $y_{2}$ is the copy of $a_{1}$ in $F_{1}$ and $\widetilde{G_{2}}$, 
       respectively. Then $(\widetilde{G_{1}})^{2}$ contains a hamiltonian cycle 
       $\widetilde{C_{1}}$ with the required properties by Theorem 
       \ref{ConnectTheorem} Ib) (using $C_{1a}$ and $\widetilde{C_{2}}$).

     Then $G = \widetilde{G_{1}}[z_{1}=z_{2}]F$, where $z_{1}$ and $z_{2}$ is the copy
     of $a_{1}$ in $\widetilde{G_{1}}$ and $F$, respectively, and $G^{2}$ contains 
     a hamiltonian cycle $C$ such that if $t_{G}(a_{i}) = 0$, then both edges of~$C$ 
     incident with $a_{i}$ are in $G$, if $t_{G}(a_{i}) = 1$, then exactly one edge of 
     $C$ incident with $a_{i}$ is in $G$, if $t_{G}(a_{i}) = 2$, then no edge of $C$ 
     incident with $a_{i}$ is in $G$, $i = 1, 2, ..., k$, by Theorem
     \ref{ConnectTheorem} Id) (using $\widetilde{C_{1}}$ and $C_{1b}$).
   \end{mathitem}

   \bigskip

   \emph{Case 2: $t_{G_{2}}(a_{2}) = 2$.}

    \noindent
    Let $B_{2}$ be an acyclic non-end block of $G_{2}$ containing the vertex $a_{2}$, 
    let $S$ be the subgraph of $G_{2}$ corresponding to the maximal connected subgraph 
    of $\mbox{Bl}(G_{2})$ containing the vertex corresponding to $B_{2}$ and not 
    containing the vertex $b_{2}$.

    Let $S_{1} = G_{2} - V(S - a_{2})$. By the induction hypothesis, $G_{2}$ contains 
    a~hamiltonian cycle $C_{2}$ such that no edge of $C_{2}$ incident with $a_{2}$ is 
    in $G$ and we can divide $C_{2}$ into hamiltonian cycles $C_{2a}$ in $S_{1}$ and 
    $C_{2b}$ in $S$ such that exactly one edge of $C_{2a}$ incident with $a_{2}$ is in 
    $S_{1}$ and exactly one edge of $C_{2b}$ incident with $a_{2}$ is in $S$.

    Now necessarily both edges of $C_{H}$ incident with $a_{2}$ are in $H$ (otherwise 
    $t_{G}(a_{2}) = 3$, a~contradiction). Set 
    $\widetilde{S_{1}} = S_{1}[x_{1}=x_{2}]H$, where $x_{1}$ and $x_{2}$ is the copy
    of $a_{2}$ in $S_{1}$ and $H$, respectively. Then $(\widetilde{S_{1}})^{2}$ 
    contains a hamiltonian cycle $C'$ with the required properties by Theorem 
    \ref{ConnectTheorem} Ib) (using $C_{2a}$ and $C_{H}$).

    Then $\widetilde{G_{2}} = \widetilde{S_{1}}[u_{1}=u_{2}]S$, where $u_{1}$ and 
    $u_{2}$ is the copy of $a_{2}$ in $\widetilde{S_{1}}$ and $S$, respectively, and 
    $(\widetilde{G_{2}})^{2}$ contains a hamiltonian cycle $\widetilde{C_{2}}$ with 
    the required properties  by Theorem \ref{ConnectTheorem}~Id) (using $C'$ and 
    $C_{2b}$). Moreover it is obvious that if $a_{1}$ is contained in a cyclic block
    of $\widetilde{G_{2}}$, then both edges of $\widetilde{C_{2}}$ incident with 
    $a_{1}$ are in $\widetilde{G_{2}}$, and if $a_{1}$ is contained in an acyclic 
    block of $\widetilde{G_{2}}$, then exactly one edge of $\widetilde{C_{2}}$ 
    incident with $a_{1}$ is in $\widetilde{G_{2}}$. Then we continue similarly as in 
    Subcase 1a) or 1b).
\end{Proof}

It is obvious that the conditions in Theorem \ref{MainTheorem} can be verified in 
polynomial time. From this it follows that the decision problem, if the square of 
a graph is hamiltonian, which is NP-complete in general (\cite{Nezn}), can be decided 
in polynomial time in the class of the graphs $G$ such that every vertex 
$x \in V_{\geq 3}(\mbox{Bl}(G))$ corresponds to a cut vertex of $G$, and for any two 
vertices $x, y \in V_{\geq 3}(\mbox{Bl}(G))$ it holds that
$\mbox{dist}_{\mbox{Bl}(G)}(x, y) \geq 4$.

\medskip

The following theorems are immediate corollaries of Theorem \ref{MainTheorem}.

\begin{Corollary}
 \label{StarTheorem}
 Let $G$ be a connected graph such that its block graph $\mbox{Bl}(G)$ is 
 homeo\-morphic to a star in which the center corresponds to a cut vertex $a$ of $G$. 
 Then the graph~$G^{2}$ is hamiltonian if and only if the vertex $a$ is contained in 
 at most two acyclic non-end blocks of $G$.
\end{Corollary}

\begin{Corollary}
 \label{StarTheorem2}
 If the block graph of $G$ with at least three vertices is a~star, then $G^{2}$ is 
 hamiltonian.
\end{Corollary}

Note that the graph in Figure \ref{Vzor} satisfies the assumptions of Corollary 
\ref{StarTheorem}. Therefore Corollary \ref{StarTheorem} (hence also Theorem 
\ref{MainTheorem}) does not follow from Theorem \ref{AbdFlaRyj}.

\section{A star in which the center corresponding to a block}
 Let $G$ be a connected graph such that its block graph $\mbox{Bl}(G)$ is homeomorphic
 to a star in which the center corresponds to a block $B_{c}$ of $G$. If $B_{c}$ is 
 acyclic, then $\mbox{Bl}(G)$ is a path and $G^{2}$ is hamiltonian by Theorem 
 \ref{Thomassen}.

 Let $B_{c}$ be cyclic. Let $k$ denote the number of cut vertices of $G$ in 
 $V(B_{c})$, let $v_{i} \in V(B_{c})$ be all cut vertices of $G$ in $B_{c}$, 
 $i = 1, 2, ..., k$. Let $C$ be a hamiltonian cycle in $(B_{c})^{2}$. We say that $C$ 
 is \emph{acceptable in $G$} if there are pairwise distinct edges 
 $v_{i}w_{i} \in E(C)$ such that $v_{i}w_{i} \in E(B_{c})$, for any 
 $i = 1, 2, ..., k$.

The following theorem gives only a sufficient condition for hamiltonicity in this 
class of graphs (in comparison with Theorem \ref{MainTheorem}).

 \begin{Theorem}
  \label{MainTheorem2}
  Let $G$ be a connected graph such that its block graph $\mbox{Bl}(G)$ is
  homeo\-morphic to a star in which the center corresponds to a block $B_{c}$ of $G$. 
  If $(B_{c})^{2}$ contains an acceptable cycle in $G$, then $G^{2}$ is hamiltonian.
 \end{Theorem}

 \begin{Proof}
  Let $v_{i} \in V(B_{c})$ be all cut vertices of $G$ in $V(B_{c})$,  
  $i = 1, 2, ..., k$. We prove the following slight strengthening of Theorem 
  \ref{MainTheorem2}.

  \medskip

  \emph{Let $G$ be a connected graph such that its block graph $\mbox{Bl}(G)$ is
  homeomorphic to a star in which the center corresponds to a block $B_{c}$ of $G$. If 
  $(B_{c})^{2}$ contains an acceptable cycle~$C$ in $G$, then $G^{2}$ contains 
  a hamiltonian cycle containing all edges of $C$ except the edges 
  $v_{i}w_{i}$, $i = 1, 2, ..., k$.}

  \medskip

  We prove this assertion by induction on $k$.

  \noindent
  \emph{1.} Let $k = 0$. Then $G = B_{c}$ and the assertion is true.

  \noindent
  \emph{2.} Suppose the assertion is true for each graph such that the block $B_{c}$ 
  contains at most $k - 1$ cut vertices, let $G$ be a graph such that its block graph 
  is homeomorphic to a star in which the center corresponds to a block $B_{c}$ of $G$ 
  and $B_{c}$ contains $k$ cut vertices of $G$, $k \geq 1$.

  Let $d \in V(\mbox{Bl}(G))$ be the vertex corresponding to $B_{c}$, let $G'$ be the 
  subgraph of $G$ such that $\mbox{Bl}(G') = \mbox{Bl}(G) - d$ and let $H_{i}$ be the 
  component of $G'$ such that $v_{i} \in V(H_{i})$, for $i = 1, 2, ..., k$. The block 
  graph $\mbox{Bl}(H_{i})$ is a path and therefore either $(H_{i})^{2}$ is hamiltonian 
  and $v_{i} \in V_{[1]}(H_{i}) \bigcup V_{[2]}(H_{i})$ (either by Theorem 
  \ref{ZobThomassen} or by Theorem \ref{Fleischner}) or $H_{i}$ is isomorphic to one 
  edge.

  Let $G_{1} = G - V(H_{1} - v_{1})$ and let $x_{1}$ and $x_{2}$ denote the copy of 
  $v_{1}$ in $G_{1}$ and $H_{1}$, respectively. Then $G = G_{1}[x_{1}=x_{2}]H_{1}$. 
  Let $C$ be an acceptable cycle in $G$. Then $C$ is also acceptable in $G_{1}$ and 
  therefore $(G_{1})^{2}$ contains a hamiltonian cycle $C_{1}$ containing all
  edges of $C$ except $v_{i}w_{i}$, $i = 2, 3, ..., k$, by the induction hypothesis. 
  The cycle $C$ is acceptable in~$G$ and therefore $x_{1}w_{1} \in E(G_{1})$ and 
  $x_{1}w_{1} \in E(C_{1})$. Then $x_{1} \in V_{[1]}(G_{1}) \bigcup V_{[2]}(G_{1})$.

  \emph{Case 1: $x_{2} \in V_{[1]}(H_{1}) \bigcup V_{[2]}(H_{1})$.}

  \noindent
  Then $G^{2}$ contains a hamiltonian cycle containing all edges of $C$ except 
  $v_{i}w_{i}$, $i = 1, 2, ..., k$, by Theorem \ref{ConnectTheorem} I).

  \emph{Case 2: $H_{1}$ is isomorphic to one edge.}

  \noindent
  Then $G^{2}$ contains a hamiltonian cycle containing all edges of $C$ except 
  $v_{i}w_{i}$, $i = 1, 2, ..., k$, by Theorem \ref{ConnectTheorem} either IIa) or 
  IIb).
 \end{Proof}

 \bigskip

 The following theorem is an immediate corollary of Theorem \ref{MainTheorem2}.

\begin{Corollary}
 \label{HamCyklus}
 Let $G$ be a connected graph such that its block graph $\mbox{Bl}(G)$ is
 homeomorphic to a star in which the center corresponds to a block $B_{c}$ of $G$. If 
 $B_{c}$ is hamiltonian, then $G^{2}$ is hamiltonian.
\end{Corollary}

Let us mention the following theorem by Schaar \cite{Schaar} that motivates the next 
conjecture.

\begin{Theoremcite} {Schaar}
 \label{Schaar}
  For every block $G$ with $|V(G)| \geq 4$ there exists a~hamiltonian cycle in $G^{2}$
  containing at least four edges of $G$.
 \end{Theoremcite}

\begin{Conjecture}
 \label{Hypoteza}
 Let $G$ be a connected graph such that its block graph $\mbox{Bl}(G)$ is
 homeo\-morphic to a star in which the center $c$ corresponds to a block $B_{c}$ of 
 $G$. If $d_{\mbox{Bl}(G)}(c) \leq~k$, $k < 5$, then $G^{2}$ is hamiltonian.
\end{Conjecture}

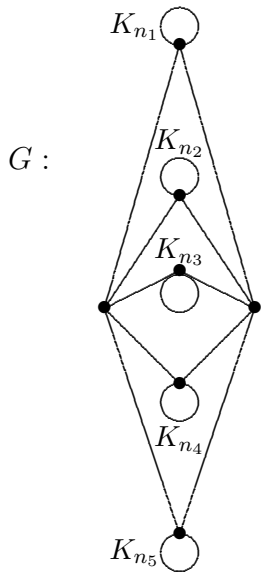
 \begin{figure}[ht]
  $$\beginpicture
   \setcoordinatesystem units <1mm,1mm>
   \setplotarea x from -30 to 30, y from -20 to 20
    \put{$G:$} at -20 20
    \put{$\bullet$} at -10   0
    \put{$\bullet$} at  10   0
    \put{$\bullet$} at   0   5
    \put{$\bullet$} at   0  15
    \put{$\bullet$} at   0 -10
    \put{$\bullet$} at   0 -30
    \put{$\bullet$} at   0  35
    \plot -10 0 0 15 10 0 0 5 -10 0 0 -10 10 0  0 -30 -10 0 0 35 10 0 /
    \circulararc 360 degrees from   0 15 center at  0 17.5
    \circulararc 360 degrees from  0 -10 center at  0 -12.5
    \circulararc 360 degrees from   0 4.5 center at   0   2
    \circulararc 360 degrees from   0 -30 center at   0  -32.5
    \circulararc 360 degrees from   0  35 center at   0   37.5
    \put{\small{$K_{n_{1}}$}} at  -6  37.5
    \put{\small{$K_{n_{2}}$}} at   0  22
    \put{\small{$K_{n_{4}}$}} at   0 -17.5
    \put{\small{$K_{n_{3}}$}} at   0   7.5
    \put{\small{$K_{n_{5}}$}} at  -6 -32.5
  \endpicture$$
  \caption{}
  \label{Priklad}
 \end{figure}

 Conjecture \ref{Hypoteza} is true for $k \leq 2$ but for $k \in \{3, 4\}$ is an open 
 problem. It is not possible to specify four edges from Theorem \ref{Schaar} and 
 therefore Conjecture~\ref{Hypoteza} is not an immediate corollary of Theorem 
 \ref{Schaar}. If Conjecture \ref{Hypoteza} is true, then the upper bound is sharp 
 as can be seen from Figure \ref{Priklad}.

\section{Conclusion and acknowledgement}
In conclusion I would like to thank Z. Ryj\'{a}\v{c}ek for valuable notes and 
recommendations.

\end{document}